\documentclass[11pt]{amsart}
\usepackage{amsmath,amsfonts,amsthm,amssymb}

\newtheorem{thm}{Theorem}
\newtheorem{lem}[thm]{Lemma}
\newtheorem{prop}[thm]{Proposition}
\newtheorem{corol}[thm]{Corollary}

\renewcommand{\i}{\infty}
\renewcommand{\a}{\alpha}
\renewcommand{\b}{\beta}
\renewcommand{\r}{\gamma}

\renewcommand{\t}{\tau}
\newcommand{\df}{\dfrac}
\newcommand{\rt}{\right}
\newcommand{\lt}{\left}
\numberwithin{equation}{section}
\numberwithin{thm}{section}
\begin{document}

\title[Mock Jacobi forms in basic hypergeometric series] {Mock Jacobi forms in \\basic hypergeometric series }
\author{Soon-Yi Kang}

\date{\today}

\address{Korea Advanced Institute for Science and Technology, Daejeon 305-701, KOREA}
\email{s2kang@kaist.ac.kr}

\subjclass[2000]{11F37, 11F50, 05A17, 33D15}
\keywords{basic hypergeometric series, mock Jacobi forms, mock modular forms, weakly holomorphic modular forms, ranks and cranks of partitions}

\begin{abstract}
We show that some $q$-series such as universal mock theta functions are linear sums of theta quotients and mock Jacobi forms of weight $1/2$, which become holomorphic parts of real analytic modular forms when they are restricted to torsion points and multiplied by suitable powers of $q$.  And we prove that certain linear sums of
$q$-series are weakly holomorphic modular forms of weight $1/2$ due to annihilation of mock Jacobi forms or completion by mock Jacobi forms.  As an application, we obtain a relation between the rank and crank of a partition.
\end{abstract}
\maketitle

\section{Introduction and summary of results}\label{intro}
A Jacobi form is a holomorphic function of two variables $z\in \mathbb{C}$ and $\t\in \mathbb{H},$ the upper half of the complex plane, which satisfies elliptic transformation properties with respect to $z$, modular transformation properties with respect to $\t$, and a certain growth condition. It is well known \cite[Theorem 1.3]{EZ} that a Jacobi form at torsion points, $z=r\t+s$, with $r$, $s\in \mathbb{Q}$, multiplied by a suitable rational power of $q=e^{2\pi i \t}$ is a modular form with respect to $\t$.  For example, the classical theta series
\begin{equation}\label{theta}\vartheta(z):=\vartheta(z;\t):=\sum_{n=-\i}^\i(-1)^{n} e^{\pi i (2n+1)z} q^{\frac{n(n+1)}{2}+\frac{1}{8}}\end{equation}
is the most famous Jacobi form and $q^{r^2/2}\vartheta(r\t+s;\t)$ is an ordinary modular form of weight $1/2$.

In his ground-breaking work on mock theta functions \cite{Zwegers}, Zwegers constructed a non-holomorphic but real analytic Jacobi form whose holomorphic part is a normalized Lerch sum and whose non-holomorphic part can be written in terms of a period integral of a weight $3/2$ unary theta series.  The holomorphic part of Zwegers' real analytic Jacobi form is defined by
\begin{equation}\label{mu}\mu(u,v):=\mu(u,v;\t):=\frac{a^{1/2}}{\vartheta(v)}\sum_{n=-\i}^\i\frac{(-b)^nq^{n(n+1)/2}}{1-aq^n},\end{equation}
where $a=e^{2\pi i u}$ and $b=e^{2\pi i v}$ for complex values $u$ and $v$.
Following Zagier, we call $\mu$ a \textit{mock Jacobi form}, because it behaves almost like a Jacobi form with two elliptic variables and when the elliptic variables are restricted to torsion points, it becomes a \textit{a mock modular form} up to a factor of a rational power of $q$.  Each mock modular form of weight $k$ is the holomorphic part of a real analytic modular form of weight $k$, of which the non-holomorphic part is determined by an ordinary modular form of weight $2-k$, called the \textit{shadow} of the mock modular form.  The real analytic modular form that has a mock modular form as its holomorphic projection is now called a harmonic weak Maass form.  And a mock theta function is a $q$-series that is a mock modular form of weight $1/2$ multiplied by a rational power of $q$ and the shadow is a unary theta series of weight $3/2$.  For a full description of the notions of mock modular forms and their shadows, the reader is referred to \cite{Zagier}.

We recall that the \textit{rank} of a partition is defined as its largest part minus its number of parts.  Then
the generating function for the number of partitions of given size and rank is given by
 $$\mathcal{R}(w;q):=\sum_{n=0}^\i \frac{q^{n^2}}{\prod_{k=1}^n(1-wq^k)(1-q^k/w)}.$$   When $w = -1$, this reduces to one of
Ramanujan's 3rd order mock theta functions $f(q)$, and Zwegers \cite{Zwegers1} proved that $q^{-1/24}f(q)$ is a component of a vector-valued real analytic modular form of weight $1/2$. With the aid of this work, Bringmann and Ono \cite{BO1,BO2} showed that $q^{-1/24}f(q)$ is the holomorphic part of a weak Maass form of weight $1/2$.  Furthermore, they \cite[Theorem 1.2]{BO2} generalized this to other roots of unity. By using Zwegers' mock Jacobi form $\mu$, Zagier \cite[Theorem 7.1]{Zagier} then found a result that subsumes and is also much simpler in formulation and proof than that by Bringmann and Ono,  as remarked by Zagier himself.  Once we represent $\mathcal{R}(w;q)$ in terms of a generalized Lambert series, we can easily deduce that
 \begin{equation}\label{CR}
\frac{q^{-1/24}\mathcal{R}(e^{2\pi i\a};q)}{e^{-\pi i \a}-e^{\pi i \a}}=\frac{\eta^3(3\t)}{\eta(\t)\vartheta(3\a;3\t)}-q^{-1/6}e^{-2\pi i \a}\mu(3\a,-\t;3\t)+q^{-1/6}e^{2\pi i \a}\mu(3\a,\t;3\t),\end{equation}
where $\eta(\t)$ is the Dedekind eta-function satisfying $q^{-1/24}\eta(\t)=\prod_{n=1}^\i(1-q^n)$.  ($q^{-1/6}$ in the last two terms of the right hand side are missing in \cite{Zagier}.) Applying Zwegers' result to (\ref{CR}), Zagier \cite[Theorem 7.1]{Zagier} proved that, for any root of unity $\zeta=e^{2\pi i\a}\neq 1$, $q^{-1/24}\mathcal{R}(\zeta;q)$ is a mock modular form of weight $1/2$ with shadow proportional to $(\zeta^{-1/2}-\zeta^{1/2})\sum_{n\in \mathbb{Z}}\lt(\frac{12}{n}\rt)n\sin(\pi n\a)q^{n^2/24}$.

In this note, we find more $q$-series which are composed of sums of theta quotients and mock Jacobi forms like the rank generating function $\mathcal{R}(w;q)$.  They become mock modular forms when the elliptic variables are restricted to torsion points. In particular, both of two \textit{universal mock-theta functions} are such series:
\begin{equation}\label{mockeven}g_2(w;q):=\sum_{n=0}^\i\frac{
(-q)_n q^{n(n+1)/2}}{(w;q)_{n+1}(q/w;q)_{n+1}},\end{equation}

\begin{equation}\label{mockodd}g_3(w;q):=\sum_{n=0}^\i\frac{
q^{n(n+1)}}{(w;q)_{n+1}(q/w;q)_{n+1}},\end{equation}
 where $(a)_0:=(a;q)_0:=1,\ \ (a)_n:=(a;q)_n:=\prod_{j=0}^{n-1}(1-aq^j),
\ \ n\geq1,$ and $(a)_\i:=(a;q)_\i:=\prod_{j=0}^{\i}(1-aq^j).$
 We call these functions universal mock-theta functions as all mock-theta functions of even (resp. odd) order are expressed in terms of $g_2$ (resp. $g_3$). This was first noticed by D. Hickerson \cite{H, H7} when he proved that mock theta conjectures on mock theta functions of orders 5 and 7 can be stated in terms of $g_3$.  Recently Gordon and McIntosh \cite{GM} observed that this is a uniform phenomenon for mock theta functions of all orders.  They noticed that $g_3$ can be written in terms of $g_2$, and hence named $g_2$ a universal mock theta function \cite{GM}.  In addition to those reviewed above, there are many results on these two functions.  However, we show that Zwegers' mock Jacobi form $\mu$ naturally reveals the properties of $g_2$ and $g_3$.
\begin{thm}\label{umt} Let $\a$ be any complex number such that $\a\neq 2n\t+m$ for $m$, $n\in Z$. If $w=e^{2\pi i \a}$, then
\begin{equation}\label{memjf}wg_2(w;q)=\frac{\eta^4(2\t)}{\eta^2(\t)\vartheta(2\a;2\t)}+wq^{-1/4}{\mu}(2\a,\t;2\t),\end{equation}
which is a mock modular form of weight $1/2$ times a rational power of $q$ at torsion points.  In particular, if $\zeta\neq 1$ is a root of unity, then $\zeta g_2(\zeta;q)$ is a mock modular form of weight $1/2$ with shadow proportional to $\sum_{n=-\i}^\i(-1)^n n \zeta^{-2n}q^{n^2}.$
\end{thm}

When $w=q^r$ for a rational number $r$,  Bringmann, Ono, and Rhoades \cite[Theorem 1.3]{BOR} proved that $\displaystyle{q^{r(1-r)}(g_2(w;q)+g_2(-w;q))}$ is a weight $1/2$ weakly holomorphic modular form on a subgroup of the full modular group $SL_2(\mathbb{Z})$ depending on $r$ by showing that the summands
are distinct Maass forms with equal non-holomorphic parts of opposite signs.  If we use elliptic transformation properties of $\mu$ and $\vartheta$, we can easily deduce a more general result from (\ref{memjf}):
 \begin{thm}\label{Anng2} Let $\a$ be any complex number such that $\a\neq 2n\t+m$ for $m$, $n\in Z$. If $w=e^{2\pi i \a}$, then we have a meromorphic Jacobi form
  \begin{equation}\label{ann}w(g_2(w;q)+g_2(-w;q))=\frac{2\eta^4(2\t)}{\eta^2(\t)\vartheta(2\a;2\t)}.\end{equation}
If $\a=r\t+s$ with $r$, $s\in \mathbb{Q}$, then it is a weight $1/2$ weakly holomorphic modular form on a subgroup of $SL_2(\mathbb{Z})$ depending on $r$ when it is multiplied by $q^{-r^2}$. \end{thm}
By writing $g_3(w;q)$ as a sum of a theta quotient and mock Jacobi forms, we can derive an analogous result for $g_3(w;q)$ as well:
\begin{thm}\label{Anng3} Let $\a$ be any complex number such that $\a\neq 3n\t+m$ for $m$, $n\in Z$. If $w=e^{2\pi i \a}$ and $\zeta_3=e^{2\pi i/3}$, then we have a meromorphic Jacobi form
  \begin{equation}\label{anng3}w^{3/2}q^{-1/24}(g_3(w;q)+g_3(\zeta_3w;q)+g_3(\zeta_3^2w;q))=\frac{3\eta^3(3\t)}{\eta(\t)\vartheta(3\a;3\t)}.\end{equation}
  If $\a=r\t+s$ with $r$, $s\in \mathbb{Q}$, then it is a weight $1/2$ weakly holomorphic modular form on a subgroup of $SL_2(\mathbb{Z})$ depending on $r$ when it is multiplied by  $q^{-3r^2/2}$.\end{thm}
We can apply Theorem \ref{Anng3} to the theory of partitions.  Let $p(n)$ denote the number of partitions of $n$.  The rank of a partition gives combinatorial interpretations for Ramanujan's famous congruences $p(5n+4)\equiv 0\pmod 5$ and
$p(7n+5)\equiv 0\pmod 7$ by splitting the partitions of $5n+4$
(resp. $7n+5$) into $5$ (resp. $7$) equal classes when sorted
according to the rank modulo $5$ (resp. $7$).  Besides the rank, there is another partition statistics, called the \textit{crank}, which gives combinatorial interpretations not only for the two congruences above but also for Ramanujan's
third partition congruence $p(11n+6)\equiv 0\pmod {11}$.   From Theorem \ref{Anng3}, we derive a relation between ranks and cranks which yields that
\begin{thm}\label{Rank-Crank} Let $N(m,n)$ and $M(m,n)$ denote the numbers of partitions of a nonnegative integer $n$ with rank $m$ and crank $m$, respectively, and $p_3(n)$ be the number of partitions of $n$ into parts congruent to $\pm 1$ modulo $3$. Then for each integer $m$, we have
\begin{equation}N(3m-1,n)+N(3m,n)+N(3m+1,n)=\sum_{0\leq 3k\leq n} M(m,k)p_3(n-3k).\end{equation}
\end{thm}

By the way, Theorem \ref{Anng2} is an immediate consequence of Ramanujan's $_1\psi_1$ summation formula \cite[p.~502]{AAR},
\begin{equation}\label{1psi1}
_1\psi_1(\a;\b;z):=\sum_{n=-\i}^\i\frac{(\a)_n}{(\b)_n}z^n
=\frac{(\b/\a)_\i(\a z)_\i(q/(\a z))_\i(q)_\i}{(q/\a)_\i(\b/(\a
z))_\i(\b)_\i(z)_\i},  \ \mathrm{for}\ |\b/\a|<|z|<1.\end{equation}
After some manipulations, we can observe that $g_2(w;q)$ arises from a positive indexed summation while $g_2(-w;q)$ arises from a negative indexed summation of the the bilateral series in (\ref{1psi1}).  This tells us the beauty of $_1\psi_1$ summation formula comes from annihilation of mock Jacobi forms in some specific cases.  The following theorem, however, shows another specific case of (\ref{1psi1}) holds due to a different reason:
\begin{thm}\label{q2cq} Define
$$K'(w;\t):=\sum_{n=0}^\i\frac{(-1)^nq^{n^2}(q;q^2)_n}{(wq^{2};q^2)_{n}(q^{2}/w;q^2)_{n}},$$
$$K''(w;\t):=\sum_{n=1}^\i\frac{(-1)^n q^{n^2}(q;q^2)_{n-1}}{(wq;q^2)_{n}(q/w;q^2)_{n}}.$$
Let $\a$ be any complex number such that $\a\neq m\t+n$ for $m$, $n\in Z$. If $w=e^{2\pi i \a}$, then we have a meromorphic Jacobi form
\begin{equation}\label{q2cqe}\frac{q^{-1/8}}{(w)^{-1/2}-(w)^{1/2}}K'(w;\t)+q^{-1/8}(w^{1/2}-w^{-1/2})K''(w;\t)=\frac{\eta^4(\t)}{\eta^2(2\t)\vartheta(\a;\t)}.\end{equation}
If $\a=r\t+s$ with $r$, $s\in \mathbb{Q}$, then it is a weight $1/2$ weakly holomorphic modular form on a subgroup of $SL_2(\mathbb{Z})$ depending on $r$ when it is multiplied by  $q^{-r^2/2}$.\end{thm}
The two terms in the left hand side of (\ref{q2cqe}) are mock Jacobi forms and each of them is about the half of a meromorphic Jacobi form in the right hand side. That is, there occurs a completion to a meromorphic Jacobi form by mock Jacobi forms not annihilation of mock Jacobi forms as in Theorem \ref{umt}.

This paper is presented as follows.  In Section 2, we review Zwegers' results on mock Jacobi forms and Zagier's notion on shadows of mock modular forms very briefly.  In Section 3, we decompose universal mock theta functions as sums of theta quotients and mock Jacobi forms and show they are mock modular forms at torsion points when they are multiplied by suitable powers of $q$. And using the results, we prove Theorems \ref{Anng2} and \ref{Anng3}.  In Section 4, using Theorem \ref{Anng3}, we establish a relation between rank and crank generating functions and discuss further relations between ranks and cranks, including Theorem \ref{Rank-Crank}.  Finally, in Section 5, we reprove some of the results in Section 3 using the theory of basic hypergeometric series and prove Theorem \ref{q2cq}.

\section{Properties of Mock Jacobi Forms and Mock Modular Forms}
We first recall important properties of a mock Jacobi form in (\ref{mu}) stated in Zwegers' Ph.D Thesis \cite{Zwegers}.
Note that our $\vartheta$ is $-i$ times of Zwegers'  $\vartheta$ function. (See \cite[Proposition 1.3]{Zwegers}.) It is well known that
the theta series in (\ref{theta}) is an odd function of $z$ having a product representation
\begin{equation}\label{jtp}\vartheta(z;\t)=q^{1/8}x^{-1/2}\prod_{n=1}^\i(1-q^n)(1-xq^{n-1})(1-x^{-1}q^{n}),\end{equation}
where $x=e^{2\pi i z}$ by the Jacobi triple product identity. 
\begin{prop}\cite[Proposition 1.4 and 1.5]{Zwegers}\label{mup} The mock Jacobi form $\mu(u,v;\t)$ in (\ref{mu}) is a meromorphic function of $u$ with simple poles in the points $u=n\t+m$ for $n,$ $m\in \mathbb{Z}$, symmetric in $u$ and $v$, and satisfies elliptic transformation properties:
\begin{equation}\mu(u+1,v)=-\mu(u,v)=e^{-2\pi i(u-v)-\pi i \t}\mu(u+\t,v)+e^{-\pi i (u-v)-\pi i \t/4},\end{equation}
and modular transformation properties:
\begin{equation}e^{\pi i /4}{\mu}(u,v;\t+1)=-(-i\t)^{-1/2}e^{\pi i (u-v)^2/\t}{\mu}(\frac{u}{\t},\frac{v}{\t};\frac{-1}{\t})+\frac{1}{2}h(u-v;\t)={\mu}(u,v;\t),\end{equation}
where $\displaystyle{h(z):=h(z;\t):=\int_\mathbb{R}\frac{e^{\pi i \t x^2-2\pi z x}}{\cosh{\pi x}}\, dx,}$ a Mordell integral.
\end{prop}
Zwegers then constructed a non-holomorphic function $R(z;\t)$ for $z\in \mathbb{C}$ so that \begin{equation}\label{zw}\widehat{\mu}(u,v;\t)=\mu(u,v;\t)-\frac{1}{2}R(u-v;\t)\end{equation} is a real analytic Jacobi form.
For rational numbers $a$ and $b$, $q^\lambda\mu(a\t,b;\t)$ is a mock modular form of weight $1/2$ for a certain rational number $\lambda$ and its shadow  is a weight $3/2$ unary theta series.

In general, the shadow $g(\t)$ of a mock modular form $f(\t)$ of weight $k$ is proportional to $\displaystyle{y^{k}\frac{\overline{\partial}}{\partial \bar{\t}}}(f+E)(\t)$, where $y=\mathrm{Im}(\t)$ and $E(\t)$ is the associated non-holomorphic function to $f(\t)$ such that $f(\t)+E(\t)$ is a real analytic modular form of weight $k$.  Under this differential operator, the holomorphic part is mapped to zero, and thus the shadow is determined by the correction function $E(\t)$ as follows:
\begin{prop}\cite[Section 5]{Zagier}\label{za} If $E(\t)$ is the correction function of a mock modular form of weight $k$ with the shadow $g(\t)$ and $y=\mathrm{Im}(\t)$, then $\displaystyle{y^{k}\frac{{\partial}}{\partial \bar{\t}}}E(\t)$ is proportional to $\overline{g(\t)}$ and
\begin{equation}\label{Err} E(\t)=(i/2)^{k-1}\int_{-\bar{\t}}^{i\i}\frac{\overline{g(-\bar{z})}}{\sqrt{-i(z+\t)}}\, dz.\end{equation}
\end{prop}
\noindent The differential operator $\displaystyle{y^{k}\frac{\overline{\partial}}{\partial \bar{\t}}}$ is also discussed as a map from the space of harmonic weak Maass forms of weight $k$ to the space of weakly holomorphic modular forms of weight $2-k$ in \cite[Proposition 3.2]{BF}.
In finding the shadow of $q^\lambda\mu(a\t,b;\t)$ then, the following result of Zwegers would be useful.
\begin{prop}\cite[(1.5)]{Zwegers}\label{dif} If $a$, $b\in \mathbb{R}$, then
 \begin{equation}\label{difr}\frac{\partial}{\partial \bar{\t}}R(a\t-b;\t)=-\frac{i}{\sqrt{2y}}e^{-2\pi a^2 y}\sum_{\nu\in \frac{1}{2}+\mathbb{Z}}(-1)^{\nu-1/2}(\nu+a)e^{-\pi i \nu^2 \bar{\t}-2\pi i \nu (a\bar{\t}-b)}.\end{equation}
\end{prop}
But if $a\in (-\frac{1}{2}, \frac{1}{2})$ and $b\in \mathbb{R}$, we may use another formula for  $R(a\t-b)$ that writes $R(a\t-b)$ in terms of the period integral of a unary theta function of weight $3/2$.
\begin{prop}\cite[Theorem 1.16]{Zwegers}\label{Rcor}  For $a$, $b\in \mathbb{R}$, define a unary theta series of weight $3/2$
\begin{equation}\label{theta32}g_{a,b}(\t):=\sum_{\nu\in a+\mathbb{Z}}\nu e^{\pi i \nu^2\t+2\pi i \nu b}.\end{equation}
Then, for $a\in (-\frac{1}{2}, \frac{1}{2})$ and $b\in \mathbb{R}$,
\begin{equation}\label{cor}\int_{-\bar{\t}}^{i\i}\frac{g_{a+\frac{1}{2},b+\frac{1}{2}}(z)}{\sqrt{-i(z+\t)}}\, dz=-e^{-\pi i a^2\t+2\pi i a(b+\frac{1}{2})}R(a\t-b).
\end{equation}\end{prop}
Although (\ref{cor}) holds for $a\in (-\frac{1}{2},\frac{1}{2})$ only, we can extend the result to all real numbers $a$ except $a\in -1/2+\mathbb{Z}$ thanks to the periodicity of $g$ as a function of $a$:
\begin{prop}\cite[Proposition 1.15]{Zwegers}\label{gp}
$g_{a,b}$ satisfies
\begin{equation}g_{a+1,b}(\t)=g_{a,b}(\t)=e^{-2\pi i a}g_{a,b+1}(\t)=-g_{-a,-b}(\t).\end{equation}
\end{prop}

\section{Universal mock theta functions in terms of mock Jacobi forms}

We first recall generalized Lambert series expansions of universal mock theta functions.  As the standard method of deriving a generalized Lambert series expansion of a $q$-series is using a limiting case of the Watson-Whipple transformation formula, we state the formula here.  For convenience, we use the notation
$(a_1,a_2,\cdots,a_k)_n:=(a_1)_n(a_2)_n\cdots(a_k)_n.$
Letting $N\to \i$ in \cite[eq.~(2.5.1),~p.~43]{GR}, we have
\begin{eqnarray}\label{watson} \sum_{n=0}^\i\frac{(\a,\b,\r,\delta,\epsilon)_n(1-\a q^{2n})q^{n(n+3)/2}}
{(\a q/\b,\a q/\r,\a q/\delta ,\a q/\epsilon,q)_n(1-\a)} \lt(-\df{\a^2}{\b\r\delta\epsilon}\rt)^n \cr
=\frac{(\a q ,\a q/(\delta\epsilon))_\i}{(\a q/\delta, \a
q/\epsilon)_\i} \sum_{n=0}^\i\frac{(\delta,\epsilon, \a
q/(\b\r))_n}{(\a q/\b,\a q/\r,q)_n}\lt(\frac{\a q}{\delta\epsilon}\rt)^n.\end{eqnarray}
Setting $\a=q$, $\b\r=q^2$, $\delta=x$, and $\epsilon=q/x$ in (\ref{watson}) and then multiplying both sides of the resulting equation by $(1-q)/((1-w)(1-q/w))$ yields the following well known Lambert series expansion of a theta quotient.
\begin{lem}\label{thetaLam} Let $\a$ be any complex number such that $\a\neq n\t$ for $n\in Z$. If $w=e^{2\pi i \a}$, then
\begin{equation}\sum_{n=-\i}^\i\frac{(-1)^nq^{n(n+1)/2}}{1-wq^{n}}=\frac{(q;q)^2_\i}{(w;q)_\i(q/w;q)_\i}.\end{equation}\end{lem}
The generalized Lambert series expansions for both of $g_2$ and $g_3$ are already known in the literature as well (for example, \cite{GM}, \cite{H}).  Utilizing (\ref{watson}) again, we have
\begin{lem}\label{uniLam} If $g_2(w;q)$ and $g_3(w;q)$ are defined as in (\ref{mockeven}) and (\ref{mockodd}), respectively, then
\begin{equation}\label{g2}g_2(w;q)=\frac{(-q)_\i}{(q)_\i}\sum_{n=-\i}^\i\frac{(-1)^nq^{n(n+1)}}{1-wq^{n}},\end{equation} and
\begin{equation}\label{g3}g_3(w;q)=\frac{1}{(q)_\i}\sum_{n=-\i}^\i\frac{(-1)^nq^{3n(n+1)/2}}{1-wq^{n}}.\end{equation}
\end{lem}

Now, we can decompose both $g_2(w;q)$ and $g_3(w;q)$ as sums of theta quotients and mock Jacobi forms. First let us prove Theorem \ref{umt}.

\begin{proof}[proof of Theorem \ref{umt}]
Using the fact $\displaystyle{\frac{1}{1-x}=\frac{1+x}{1-x^2}}$, we deduce  from (\ref{g2}) that
$$g_2(w;q)=\frac{(-q)_\i}{(q)_\i}\lt(\sum_{n=-\i}^\i\frac{(-1)^nq^{n^2+n}}{1-w^2q^{2n}}+w\sum_{n=-\i}^\i\frac{(-q)^nq^{n^2+n}}{1-w^2q^{2n}}\rt).$$
Using Lemma \ref{thetaLam} and the definitions of $\vartheta(z)$  and $\mu(u,v)$  in  (\ref{theta})  and (\ref{mu}), respectively, we establish (\ref{memjf}).
Hence by the results in Section 2, we can see that $w g_2(w;q)$ is a mock modular form of weight $1/2$ at torsion points when it is multiplied by a rational power of $q$. In particular, if we restrict $\a$ to a rational number, then in (\ref{memjf}), the first term is a weakly holomorphic modular form of weight $1/2$ and the second term is a mock modular form of weight $1/2$ with the correction term $-e^{2\pi i\a}q^{-1/4}R(2\a-\t;2\t)/2$.  By setting $a=-1/2$ in (\ref{difr}), we find that  \begin{equation}\label{difr2}\frac{\partial}{\partial \bar{\t}}R(-\frac{\t}{2}-b)=-\frac{i}{\sqrt{2y}}e^{\pi i \t/4+\pi i b}\sum_{n\in \mathbb{Z}}(-1)^{n}ne^{-\pi i n^2\bar{\t}+2\pi i n b},\end{equation}
and hence it follows from Proposition \ref{za} that the shadow of $e^{-\pi i b}q^{-1/4}\mu(-\t/2,b;\t)$ is
\begin{equation}\label{sh}{\frac{\sqrt{i}}{2}\sum_{n\in \mathbb{Z}}(-1)^{n}ne^{\pi i n^2{\t}+2\pi i n b}}.\end{equation}  Therefore, substituting $b=-2\a$ and replacing $\t$ by $2\t$ in (\ref{sh}) completes the proof for $\zeta g_2(\zeta;q)$, where $\zeta=e^{2\pi i\a}\neq 1$.
\end{proof}
We have similar results for $g_3$.
\begin{thm}\label{ag} Let $\a$ be any complex number such that $\a\neq 3m\t+n$ for $m$, $n\in Z$, then
\begin{equation}\label{momjf}q^{-1/24}e^{3\pi i\a}g_3(e^{2\pi i\a};q)=\frac{\eta^3(3\t)}{\eta(\t)\vartheta(3\a;3\t)}+q^{-1/6}e^{2\pi i \a}\mu(3\a,\t;3\t)+q^{-2/3}e^{4\pi i \a}\mu(3\a,2\t;3\t),\end{equation}
which is a mock modular form of weight $1/2$ times a rational power of $q$ at torsion points.  In particular, if $\zeta\neq 1$ is a root of unity, then
 $q^{-1/24}\zeta^{3/2} g_3(\zeta;q)$ is a mock modular form of weight $1/2$ with shadow proportional to $\sum_{n=-\i}^\i(\frac{12}{n}) n \sin(\pi n \a)q^{n^2/24}.$\end{thm}
\textit{Remark.} As stated in \cite[Theorem 7.1]{Zagier},
$$\mathcal{R}(w;q)=\frac{1-w}{(q)_\i}\sum_{n=-\i}^\i\frac{(-1)^nq^{(3n^2+n)/2}}{1-wq^n}.$$
By comparing this with the Lambert series expansion of $g_3(w;q)$ in (\ref{g3}), we can easily see that $\displaystyle{g_3(w;q)=-\frac{1}{w}+\frac{1}{w(1-w)}R(w;q)}$.  Hence we may obtain the result for $g_3(\zeta;q)$ above from the fact that $q^{-1/24}\mathcal{ R}(\zeta;q)$ is a mock modular form stated in the introduction as well.
\begin{proof}
The proof is similar to that of Theorem \ref{umt}.  It follows from  (\ref{g3}) and ${\frac{1}{1-x}=\frac{1+x+x^2}{1-x^n}}$ that
$$g_3(w;q)=\frac{1}{(q)_\i}\lt(\sum_{n=-\i}^\i\frac{(-1)^nq^{3n(n+1)/2}}{1-w^3q^{3n}}+\sum_{n=-\i}^\i\frac{(-q)^nwq^{3n(n+1)/2}}{1-w^3q^{3n}}
+\sum_{n=-\i}^\i\frac{(-q^2)^nw^2q^{3n(n+1)/2}}{1-w^3q^{3n}}\rt).$$  From Lemma \ref{thetaLam} and the definitions of $\vartheta(z)$  and $\mu(u,v)$  in  (\ref{theta})  and (\ref{mu}), respectively, we deduce (\ref{momjf}), where $w=e^{2\pi i \a}$.  Hence by the results in Section 2, we can see that $q^{-1/24}e^{3\pi i\a}g_3(e^{2\pi i\a};q)$ is a mock modular form of weight $1/2$ at torsion points when it is multiplied by a rational power of $q$. In particular, if we restrict $\a$ to a rational number, the first term in (\ref{momjf}) is a weakly holomorphic modular form of weight $1/2$ and the last two terms are mock modular forms of weight $1/2$ with correction terms $-\frac{e^{2\pi i\a} q^{-1/6}}{2}R(3\a-\t;3\t)$ and $-\frac{e^{4\pi i\a} q^{-2/3}}{2}R(3\a-2\t;3\t)$ by (\ref{zw}). Using Propositions \ref{Rcor} and \ref{gp}, we can find that
\begin{eqnarray}&&\frac{e^{2\pi i\a} q^{-1/6}}{2}R(3\a-\t;3\t)+\frac{e^{4\pi i\a} q^{-2/3}}{2}R(3\a-2\t;3\t)\cr
&=&-\frac{1}{2}\int_{-3\bar{\t}}^{i\i}\frac{e^{\pi i/3}g_{\frac{1}{6},-3\a+\frac{1}{2}}(z)-e^{2\pi i/3}g_{\frac{1}{6},3\a-\frac{1}{2}}(z)}{\sqrt{-i(z+3\t)}}\, dz,\end{eqnarray} and the shadow is $$\frac{1}{2}\lt(e^{\pi i/3}g_{\frac{1}{6},-3\a+\frac{1}{2}}(3\t)-e^{2\pi i/3}g_{\frac{1}{6},3\a-\frac{1}{2}}(3\t)\rt)=\frac{1}{6}\sum_{n=-\i}^\i (-1)^n (6n+1)q^{(6n+1)^2/24}\sin(\pi(6n+1)\a).$$
\end{proof}

Next, let us prove Theorem \ref{Anng2} and Theorem \ref{Anng3}.

\begin{proof}[Proof of Theorem \ref{Anng2}]
The result follows immediately from (\ref{memjf}), as we have
$$g_2(-e^{2\pi i\a};q)=g_2(e^{2\pi i(\a+\frac{1}{2})};q)=\frac{\eta^4(2\t)e^{-2\pi i\a}}{\eta^2(\t)\vartheta(2\a;2\t)}-q^{-1/4}\mu(2\a,\t;2\t)$$
by elliptic properties of $\vartheta$ and $\mu$ such as $\vartheta(z+1)=-\vartheta(z)$ and $\mu(u+1,v)=-\mu(u,v)$.
\end{proof}

\begin{proof}[Proof of Theorem \ref{Anng3}] Applying (\ref{momjf}) three times with the elliptic properties of $\vartheta$ and $\mu$, we obtain that
\begin{eqnarray}&&q^{-1/24}(g_3(w;q)+g_3(\zeta_3w;q)+g_3(\zeta_3^2w;q))=\frac{3\eta^3(3\t)e^{-3\pi i\a}}{\eta(\t)\vartheta(3\a;3\t)}\cr
&&\ +e^{-\pi i\a}q^{-1/6}\mu(3\a,\t;3\t)(1-e^{-\pi i/3}+e^{-2\pi i/3})+e^{\pi i\a}q^{-2/3}\mu(3\a,2\t;3\t)(1-e^{\pi i/3}+e^{2\pi i/3}),\nonumber\end{eqnarray}
where the last two terms are zero.
\end{proof}

\section{Applications to partitions}

Using (\ref{CR}) or using Theorem \ref{Anng3}, we also have a similar result to Theorem \ref{Anng3} for the rank generating function $\mathcal{R}(w;q)$: If $w=e^{2\pi i\a}$ for any complex number $\a\neq 3n\t+m$ with $m$, $n\in \mathbb{Z}$ and $\zeta=e^{2\pi i/3}$, then
\begin{equation}\label{AnnCR}
w^{3/2}q^{-1/24}\lt(\frac{\mathcal{R}(w;q)}{w(1-w)}+\frac{\mathcal{R}(w\zeta;q)}{w\zeta(1-w\zeta)}+\frac{\mathcal{R}(w\zeta^2;q)}{w\zeta^2(1-w\zeta^2)}\rt)
=\frac{3\eta^3(3\t)}{\eta(\t)\vartheta(3\a;3\t)}.\end{equation}

We recall the \textit{crank} of a partition, which is defined by
$$\mathrm{crank}=\left\{%
\begin{array}{ll}
    \mathrm{the\ largest\ part}, & \hbox{{if} part $1$ does not appear;} \\
    e(k)-k, & \hbox{{if} part $1$ appears $k$ times,} \\
\end{array}%
\right.$$ where $e(k)$ denotes the number of parts in the partition that are strictly
larger than $k$. Let $M(m,n)$ denote the number of partitions of a nonnegative integer $n$ with crank $m$.  If we define $M(0,1)=-1$, $M(-1,1)=M(1,1)=1$ and $M(m,1)=0$ otherwise, the generating function for $M(m,n)$ is given by, \cite{AG},
\begin{equation}\label{gencrank}
\mathcal{C}(w;q):=1+\sum_{m=-\i}^\i\sum_{n=1}^\i M(m,n)w^mq^n=\frac{(q;q)_\i}{(wq;q)_\i(q/w;q)_\i}.\end{equation}
Then from \eqref{AnnCR}, we deduce the following relation between $\mathcal{R}(w;q)$ and $\mathcal{C}(w;q)$:
\begin{equation}\label{rank-crank}
\frac{\mathcal{R}(w;q)}{w(1-w)}+\frac{\mathcal{R}(w\zeta;q)}{w\zeta(1-w\zeta)}+\frac{\mathcal{R}(w\zeta^2;q)}{w\zeta^2(1-w\zeta^2)}
=\frac{3(q^3;q^3)_\i\mathcal{C}(w^3;q^3)}{(q;q)_\i(1-w^3)}.\end{equation}
Furthermore, we note that the left hand side of equation \eqref{rank-crank} can be written as
\begin{equation}\label{left}
L:=\frac{1}{1-w^3}\lt(\mathcal{R}(w;q)(\frac{1}{w}+1+w)+\mathcal{R}(w\zeta;q)(\frac{1}{w\zeta}+1+w\zeta)
+\mathcal{R}(w\zeta^2;q)(\frac{1}{w\zeta^2}+1+w\zeta^2)\rt).\end{equation}
If we apply
$$\mathcal{R}(w;q)=1+\sum_{m=-\i}^\i\sum_{n=1}^\i N(m,n)w^mq^n$$  into \eqref{left}, we find that
\begin{equation}
L=\frac{3}{1-w^3}\sum_{m=-\i}^\i\sum_{n=0}^\i\lt(N(3m,n)+N(3m+1,n)+N(3m-1,n)\rt)w^{3m}q^n,\end{equation}
from the fact $1+\zeta+\zeta^2=0$.
Hence, from \eqref{rank-crank}, we derive the theorem below after replacing $w^3$ by $w$:
\begin{thm}\label{rc} Let $N(m,n)$ and $M(m,n)$ denote the numbers of partitions of a nonnegative integer $n$ with rank $m$ and crank $m$, respectively. Then we have
\begin{eqnarray}\label{rcc}&&\sum_{m=-\i}^\i\sum_{n=0}^\i\lt(N(3m,n)+N(3m+1,n)+N(3m-1,n)\rt)w^{m}q^n\cr
&&\hspace{1in}=\frac{(q^3;q^3)_\i}{(q;q)_\i}\sum_{m=-\i}^\i\sum_{n=0}^\i M(m,n)w^{m}q^{3n}.\end{eqnarray}
\end{thm}

Theorem \ref{Rank-Crank} follows from Theorem \ref{rc} immediately.  For example, when $n=4$, there are 5 different partitions, and
$N(-1,4)+N(0,4)+N(1,4)=3=M(0,0)p_3(4)+M(0,1)p_3(1)$.

For another consequence of Theorem \ref{rc}, we substitute $w=-1$ in \eqref{rcc} and obtain:
\begin{corol} Let $N(t,m,n)$ and $M(t,m,n)$ denote the number of partitions of $n$ with rank and crank, respectively, congruent to $t$ modulo $m$.  Then
 \begin{equation}\label{rcc-1}\sum_{n=0}^\i N_{(6)}(n)q^n=\frac{(q^3;q^3)_\i}{(q;q)_\i}\sum_{n=0}^\i M_{(2)}(n)q^{3n}=\frac{(q^3;q^3)^2_\i}{(q;q)_\i(-q^3;q^3)^2_\i},\end{equation}
where $\displaystyle{N_{(6)}(n):=\sum_{t=0,\pm 1}N(t,6,n)-\sum_{t=3,\pm 2}N(t,6,n)}$ and $\displaystyle{M_{(2)}(n):=M(0,2,n)-M(1,2,n)}.$

And hence, \begin{equation}\label{rcc-2}N_{(6)}(n)=\sum_{0\leq 3k\leq n} M_{(2)}(k)p_3(n-3k).\end{equation}
\end{corol}
From the work of Treneer \cite{T}, we know that $N_{(6)}(n)$ has infinitely many congruences in arithmetic progressions modulo any prime coprime to $6$. Although $M_{(2)}(n)$ satisfies a Ramanujan type congruence modulo $5$, namely, $M_{(2)}(5n+4)\equiv 0\pmod 5$ which is proved in \cite{CKL}, there seems no such simple congruences for $N_{(6)}(n)$.  However, it seems that $N_{(6)}(2n)>0$ and $N_{(6)}(2n+1)<0$ for any non-negative integers $n$ such as $M_{(2)}(n)$. But note that $N_{(6)}(1)>0$ while $M_{(2)}(1)<0$.  We will leave this observation as a conjecture.

In closing this section, we should mention that Theorem \ref{Anng3} is a special case of a more general identity due to Atkin and Swinnerton-Dyer \cite[(5.1)]{AtS} who used the general identity to prove Dyson's results for the rank modulo 5 and 7.  It also leads to the rank-crank PDE, \cite[Theorem 1.1]{AtG}, which is another relation between the rank and crank.

\section{Mock Jacobi forms that arise in $_1\psi_1$ summation formula}

 It is much more convenient to discuss our subject in terms of the following equivalent form of Ramanujan's $_1\psi_1$ summation formula (\ref{1psi1}):
\begin{prop}[{Three variable Reciprocity
Theorem} \cite{K07}]\label{3var} If $c\neq -aq^{-n}$, $-bq^{-m}$ for non-negative integers $n$ and $m$, then
\begin{equation}\label{3rec}
\rho(a,b,c)-\rho(b,a,c)=\left(\frac{1}{b}-\frac{1}{a}\right)
\frac{(c)_\i(aq/b)_\i(bq/a)_\i(q)_\i}{(-c/a)_\i(-c/b)_\i(-aq)_\i(-bq)_\i},
\end{equation}
 where
\begin{equation}\label{rho3}\rho(a,b,c):=(1+\frac{1}{b})\sum_{n=0}^\i\frac{(c)_n(-1)^nq^{n(n+1)/2}(a/b)^n}{(-aq)_n(-c/b)_{n+1}}.\end{equation}\end{prop}
If we set $c=0$ in (\ref{3rec}), we obtain Ramanujan's reciprocity theorem \cite[Theorem 1.1]{K07} and identity (\ref{ann}) is
a special case of (\ref{3rec}) when $c=-q$.

\begin{proof}[Second proof of Theorem \ref{Anng2}]\label{c=-q}
Letting $b=-a=-w$ and $c=-q$ in the three variable reciprocity theorem (\ref{3rec}), we find that
$$\sum_{n=0}^\i
\frac{(-q)_n q^{n(n+1)/2}}{(-w)_{n+1}(-q/w)_{n+1}}+\sum_{n=0}^\i
\frac{(-q)_n
q^{n(n+1)/2}}{(w)_{n+1}(q/w)_{n+1}}=\frac{2(-q;q)^2_\i(q^2;q^2)_\i}{(w^2;q^2)_\i(q^2/w^2;q^2)_\i}.$$
Applying Jacobi triple product identity (\ref{jtp}), we obtain (\ref{ann}).\end{proof}

As commented in the introduction, Bringmann, Ono, and Rhoades proved Theorem \ref{Anng2} when $w=q^r$ for a rational number $r$ in \cite[Theorem 1.3]{BOR}. The theorem has one more result, and we state it here again with their notations:
\begin{prop}\cite[Theorem 1.3]{BOR}\label{k} Let $K'$ and $K''$ be defined as in Theorem \ref{q2cq}.
If $\zeta_\r:=e^{2\pi i/\r}$, $f_\r:=2\r/\gcd(\r,4)$, and
$$\tilde{K}(\a,\r;\t):=\frac{1}{4}\csc(\a\pi/\r)q^{-\frac{1}{8}}K'(\zeta_\r^\a;\t)+\sin(\a\pi/\r)q^{-\frac{1}{8}}K''(\zeta_\r^\a;\t),$$
then $\tilde{K}(\a,\r;2f_\r^2\t)$ is a weight 1/2 weakly holomorphic
modular form on $\Gamma_1(64f_\r^4)$.
\end{prop}
Theorem \ref{q2cq} is more general than Proposition \ref{k} and it is an immediate consequence of the reciprocity theorem above.
\begin{proof}[1st proof of Theorem \ref{q2cq}]
Replacing $q$ by $q^2$, and letting $a=-w$, $b=-wq$ and $c=q$ in (\ref{3rec}), and simplifying the resulting equation, we deduce that
\begin{equation}\label{ktheta}\frac{1}{1-w}K'(w;\t)+(1-\frac{1}{w})K''(w;\t)=\frac{(q;q^2)^3_\i(q^2;q^2)_\i}{(w;q)_\i(q/w;q)_\i}.\end{equation}
Using Jacobi triple product identity (\ref{jtp}), we finish this proof.\end{proof}

\begin{proof}[2nd proof of Theorem \ref{q2cq}]
We first find generalized Lambert series expansions for $K'(w;\t)$ and $K''(w;\t)$ so that we can observe what happens in terms of mock Jacobi forms.
If we let $\epsilon\to\i$ and then replace $q$ by $q^2$ in (\ref{watson}), we have
\begin{eqnarray}\label{einfty} 1+\sum_{n=1}^\i\frac{(\a q^2;q^2)_{n-1}(\b,\r,\delta ;q^2)_n(1-\a q^{4n})q^{2n(n+1)}}
{(\a q^2/\b,\a q^2/\r,\a q^2/\delta ,q^2;q^2)_n} \lt(\df{\a^2}{\b\r\delta }\rt)^n \cr
=\frac{(\a q^2;q^2)_\i}{(\a q^2/\delta;q^2 )_\i} \sum_{n=0}^\i\frac{(\delta , \a
q^2/(\b\r);q^2)_nq^{n(n+1)}}{(\a q^2/\b,\a q^2/\r,q^2;q^2)_n}\lt(-\frac{\a }{\delta }\rt)^n,\end{eqnarray}
since $\lim(\epsilon)_n\epsilon^{-n}=(-1)^nq^{n(n-1)/2}.$
Setting $\a=1$, $\delta=q$, $\b=w$, and $\r=1/w$  in (\ref{einfty}), we deduce that
\begin{equation}\label{k'}K'(w;\t)=\frac{(q;q^2)_\i}{(q^2;q^2)_\i}\lt\{1+(1-w)(1-\frac{1}{w})\sum_{n=1}^\i\frac{(1+q^{2n})q^{2n^2+n}}{(1-q^{2n}/w)(1-wq^{2n})}\rt\}.\end{equation}
Hence \begin{equation}\label{k'1}\frac{1}{1-w}K'(w;\t)=\frac{(q;q^2)_\i}{(q^2;q^2)_\i}\sum_{n=-\i}^\i\frac{q^{2n^2+n}}{1-wq^{2n}}.\end{equation}
Next, after letting $\epsilon\to\i$, setting $q$ by $q^2$, $\a=q^2$, $\delta=q$, $\b=wq$, and $\r=q/w$  in (\ref{watson}), we find that
\begin{equation}\label{k''}K''(w;\t)=-\frac{(q;q^2)_\i}{(q^2;q^2)_\i}\sum_{n=0}^\i\frac{(1+q^{2n+1})q^{2n^2+3n+1}}{(1-q^{2n+1}/w)(1-wq^{2n+1})}.\end{equation}
Thus \begin{equation}\label{k''1}(1-\frac{1}{w})K''(w;\t)=-\frac{(q;q^2)_\i}{(q^2;q^2)_\i}\sum_{n=-\i}^\i\frac{q^{2n^2+3n+1}}{1-wq^{2n+1}}.\end{equation}
Adding (\ref{k'1}) and (\ref{k''1}), we derive that
\begin{equation}\label{k'addk''}\frac{1}{1-w}K'(w;\t)+(1-\frac{1}{w})K''(w;\t)=\frac{(q;q^2)_\i}{(q^2;q^2)_\i}\sum_{n=-\i}^\i\frac{(-1)^n q^{n(n+1)/2}}{1-wq^{n}}.\end{equation}
Now by using Lemma \ref{thetaLam} and Jacobi triple product identity (\ref{jtp}), we complete the proof.\end{proof}

According to (\ref{k'1}) and (\ref{k''1}), each of $\displaystyle{\frac{1}{1-w}K'(w;\t)}$ and $\displaystyle{(1-\frac{1}{w})K''(w;\t)}$ belongs to a more general type of mock Jacobi forms than $\mu$ in (\ref{mu}), called a level 2 Appell function by Zagier and Zwegers \cite{Zwegers2}.  The proof above shows that when the index $n$ in the series in the right hand side of (\ref{k'addk''}) is even, the sum equals $\displaystyle{\frac{1}{1-w}K'(w;\t)}$ and when it is odd, the sum is $\displaystyle{(1-\frac{1}{w})K''(w;\t)}$.  Hence the theta quotient in (\ref{q2cqe}) is obtained by adding two mock Jacobi forms not by canceling two mock jacobi forms such as in Theorems \ref{Anng2} and \ref{Anng3}.

\section{Concluding Remarks}

We have observed that once we write $q$-series as generalized Lambert series

$\displaystyle{\frac{1}{\theta(\t;k\t)}\sum_{n=-\i}^\i\frac{(-1)^nq^{kn(n+1)/2}}{1-wq^{n}}}$
($k<3$), Zwegers' work on mock Jacobi forms provide short but beautiful and complete information on the $q$-series. When $k=1$ and  $\theta(\t;k\t)=1$, the series is a theta quotient and when $k=2$ or $3$ and $\theta=\vartheta$, it is a sum of a theta quotient and $k-1$ mock Jacobi forms.  For larger values of $k$, Zagier and Zwegers \cite{Zwegers2} observed that it is still a holomorphic part of a mock modular form which is a sum of a theta quotient and $k-1$ correction terms.
Further study in this direction may help us to see which $q$-series is a mock modular form.

\end{document}